\documentclass[oneside]{amsart}
\usepackage[T1]{fontenc}
\usepackage[latin1]{inputenc}
\usepackage{float}
\usepackage{amssymb}

\makeatletter

 \theoremstyle{plain}
\newtheorem{thm}{Theorem}[section]
  \theoremstyle{plain}
  \newtheorem{lem}[thm]{Lemma}
  \theoremstyle{definition}
  \newtheorem{defn}[thm]{Definition}
  \theoremstyle{remark}
  \newtheorem{rem}[thm]{Remark}
  \theoremstyle{plain}
  \newtheorem{prop}[thm]{Proposition}
  \theoremstyle{plain}
  \newtheorem{cor}[thm]{Corollary}
 \theoremstyle{definition}
  \newtheorem{example}[thm]{Example}

\usepackage{euscript}
\usepackage{eufrak}
\usepackage{dcolumn}
\newcolumntype{d}[1]{D{.}{.}{#1}} 
\newcommand{\SLZ}{SL(\ensuremath{2,\mathbb{Z}})}
\newcommand{\PSLR}{PSL(\ensuremath{2,\mathbb{R}})}
\newcommand{\PSLZ}{PSL(\ensuremath{2,\mathbb{Z}})}

\newcommand{\GLZ}{GL(\ensuremath{2,\mathbb{Z}})}

\newcommand{\sgn}{\text{sgn}}
\newcommand{\Arg}{\text{Arg}}
\renewcommand{\H}{\mathcal{H}}
\newcommand{\MAS}{\ensuremath{\EuScript{M}}}
\newcommand{\Mas}[1]{\EuScript{M}(#1)}

\renewcommand{\mod}{\  \text{mod} \ }
\newcommand{\T}[2][v]{T_{#2,k}^{#1}}
\newcommand{\Ta}[2][v]{T_{#2,k}^{#1*}}
\newcommand{\Th}[2][v]{\Theta_{#2,k}^{#1}}

\address{
Institut f\"ur Theoretische Physik \\
TU-Clausthal  \\
Abteilung Statistische Physik und Nichtlineare Dynamik \\
Arnold-Sommerfeld-Stra{\ss}e-6 \\
38678 Clausthal-Zellerfeld\\
Germany}
\email{fredrik@tu-clausthal.de}

\makeatother
\begin{document}
\title[Hecke Operators for Integer Weights]{Hecke Operators for Maass waveforms on $\PSLZ$ with Integer Weight and Eta Multiplier.}

\author{Fredrik Strömberg}

\begin{abstract}
We construct Hecke operators acting on Maass waveforms of integer
non-zero weight and transforming according to a non-trivial multiplier
system on the modular group. Using these Hecke operators we obtain
multiplicativity relations for the Fourier coefficients of such Maass
waveforms. These relations generalize the usual weight zero relations.
We also obtain an unexpected relation between coefficients with positive
and negative indices with a constant of proportionality involving
the Laplace eigenvalue. Numerical examples of multiplicativity relations
are given at the end of the paper.
\end{abstract}
\maketitle

\section{Introduction }

In the theory of modular forms Hecke operators are extremely useful
in order to obtain certain arithmetical information about the objects
of interest. The most common application is to derive multiplicativity
relations between Fourier coefficients of modular forms thus demonstrating
the existence of Euler product expansions.

Introduced by Hecke \cite{Hecke:37} and subsequently developed by
e.g.~Shimura \cite{shimura}, the theory of Hecke operators, $T_{n}$,
has been developed in great detail. Foremost this was done in the
setting of the modular group, integer weight and the multiplier being
trivial or a Dirichlet character. For congruence subgroups $\Gamma_{0}\left(m\right)$
the theory was worked out by Atkin and Lehner \cite{atkinlehner}.
In these cases Hecke operators $T_{p}$ are defined for all primes
and extended to all positive integers by multiplicativity. On $\Gamma_{0}\left(4\right)$
Shimura \cite{shimura:73:half_integral} showed that there only exist
Hecke operators $T_{p^{2}}$ for $p\ge2$, acting on modular forms
of half integer weight and transforming according to the theta multiplier
system. For more introductory texts see for example \cite{apostol,MR2112196,gunning,koblitz-introduction-to-elliptic}.

Another approach to Hecke operators acting on spaces of modular functions
with a general multiplier system and weight was initiated by Wohlfahrt
\cite{MR0106888} and subsequently developed further and applied by
van Lint \cite{van-Lint:HeckeOperators:MR0090616}. In this paper
we will see how these operators can be used to derive new multiplicativity
relations for Fourier coefficients of Maass waveforms with integer
weights and eta multiplier system.

\section{Basic Concepts and Notation}

Although most of our definitions would work also for subgroups, we
will restrict ourselves to the modular group, $\Gamma=\PSLZ\approx\left\{ \pm I\right\} \backslash\overline{\Gamma}$,
where $\overline{\Gamma}=\SLZ$ is the group of two-by-two matrices
with integer entries and determinant one and $I$ is its identity
element. The hyperbolic upper half-plane is defined as $\H=\left\{ z\in\mathbb{C}\,|\,\Im z>0\right\} $
endowed with the metric $ds=\frac{\left|dz\right|}{y}$. $\Gamma$
is known to act as the group of orientation preserving isometries
of $\H$ and the action of $\Gamma$ on functions $f\left(z\right)$
defined on $\H$ is given by the so called the slash action $g:f\left(z\right)\mapsto f_{|g}\left(z\right)=f\left(gz\right),$
where $gz=\frac{az+b}{cz+d}$ for $g=\left(\begin{smallmatrix}a & b\\
c & d\end{smallmatrix}\right)\in\Gamma$. For any $z\in\mathbb{C}$ we always consider the principal branch
of the argument, $\Arg z\in]-\pi,\pi]$ and set $e\left(z\right)=e^{2\pi iz}$.
The greatest common divisor symbol $\left(\cdot,\cdot\right)$ is
extended to $0$ by $\left(a,0\right)=a$ for all integers $a$. 

If $k$ is an integer we know from the theory of holomorphic modular
forms (e.g.~\cite{gunning} or \cite{rankin:mod}) that there is
a weight-$k$ action of $\overline{\Gamma}$ on holomorphic functions
$F\left(z\right)$ on $\H$ given by $g:F\left(z\right)\mapsto F_{|\left\{ k,g\right\} }\left(z\right)=\left(cz+d\right)^{-k}F\left(gz\right)$
for $g=\left(\begin{smallmatrix}a & b\\
c & d\end{smallmatrix}\right)\in\overline{\Gamma}$. We are interested in non-holomorphic modular forms, Maass waveforms,
which are eigenfunctions of the weight $k$ Laplacian: \[
\Delta_{k}=\Delta-iky\frac{\partial}{\partial x}=y^{2}\left(\frac{\partial^{2}}{\partial x^{2}}+\frac{\partial^{2}}{\partial y^{2}}\right)-iky\frac{\partial}{\partial x}.\]
$\Delta_{k}$ can be derived by separation of variables from the Laplace-Beltrami
operator corresponding to a certain metric on $\H\times S^{1}$, cf.
e.g.~Selberg \cite[pp.\ 81-83]{MR0088511} or Maass \cite[pp.\ 174-175]{maass:modular_functions}.
If $\chi$ is any (unitary) representation of $\overline{\Gamma}$
and $F$ is a holomorphic modular form satisfying $F_{|\left\{ k,g\right\} }\left(z\right)=\chi\left(g\right)F\left(z\right)$
for all $g\in\overline{\Gamma}$ the function $f\left(z\right)=\left(\Im z\right)^{\frac{k}{2}}F\left(z\right)$
is an eigenfunction of $\Delta_{k}$ and satisfies \[
f_{|\left[k,g\right]}\left(z\right)=\chi\left(g\right)f\left(z\right),\,\mbox{for }g=\left(\begin{smallmatrix}a & b\\
c & d\end{smallmatrix}\right)\in\overline{\Gamma},\]
 where $f_{|\left[k,g\right]}\left(z\right)=j_{g}\left(z;k\right)^{-1}f\left(gz\right),$
and $j_{g}\left(z;k\right)=e^{ik\Arg\left(cz+d\right)}.$ For any
real number $k$, a multiplier system of weight $k$ on $\overline{\Gamma}$
is a function $v:\overline{\Gamma}\rightarrow S^{1}$ such that $v\left(-I\right)=e^{-\pi ik}$
and $v\left(AB\right)=j_{A}\left(Bz;k\right)j_{B}\left(z;k\right)j_{AB}\left(z;k\right)^{-1}v\left(A\right)v\left(B\right)$
for all $A,B\in\overline{\Gamma}$. If $k$ is integer then clearly
$v\left(AB\right)=v\left(A\right)v\left(B\right)$. 

An example of a multiplier system of weight $k$ is $v_{\eta}^{2k}$,
where $v_{\eta}$ is the so called eta multiplier system given by
the Dedekind eta function $\eta\left(z\right)=e\left(\frac{z}{24}\right)\prod_{n=1}^{\infty}\left(1-e\left(nz\right)\right)$.
One can show that \[
\eta_{\left[k,g\right]}\left(z\right)=v_{\eta}\left(g\right)\,\eta\left(z\right),\, g=\left(\begin{smallmatrix}a & b\\
c & d\end{smallmatrix}\right)\in\overline{\Gamma},\]
where $v_{\eta}\left(g\right)$ is a 24-th root of unity given explicitly
either in terms of a Dedekind sum \cite[thm 3.4]{apostol} or a quadratic
residue symbol \cite[p.\ 51]{knopp:modular}. For our purposes we
use the latter approach resulting in the following formula, valid
for positive odd $c$ and integer $k$: \begin{equation}
v_{\eta}^{2k}\left(\begin{smallmatrix}\begin{array}{cc}
a & b\\
c & d\end{array}\end{smallmatrix}\right)=e\left(\frac{k}{12}\left[\left(a+d\right)c-bd\left(c^{2}-1\right)-3c\right]\right),\,\left(\begin{smallmatrix}a & b\\
c & d\end{smallmatrix}\right)\in\overline{\Gamma}.\label{eq:v_eta_explicit}\end{equation}

To distinguish between different multiplier systems, in particular
between $v$ and $\overline{v}$ the following Lemma from \cite[p.\ 11]{van-Lint:HeckeOperators:MR0090616}
is useful. 

\begin{lem}
\label{lem:mult_system_by_xi}Let $v$ be a multiplier system of weight
$k$ on $\overline{\Gamma}$. Then there exists a unique $\xi\in S^{1}$
and a function $w:\overline{\Gamma}\rightarrow\mathbb{Z}$ independent
of $k$ such that for each $A\in\overline{\Gamma}$ \[
v\left(A\right)=\xi^{w(A)}.\]
Furthermore $\xi$ is given by $\xi=v(T)$ where $T=\left(\begin{smallmatrix}1 & 1\\
0 & 1\end{smallmatrix}\right)$ and can be written as $\xi=\xi_{6}e^{\frac{\pi ik}{6}}$ for a unique
6-th root of unity, $\xi_{6}$. Thus there exist exactly 6 different
multiplier systems (of any real weight $k$) on $\Gamma$ (cf. also
\cite[p.\ 83]{rankin:mod}). 
\end{lem}
\begin{defn}
\label{def:symbol_chi}We introduce the notation $\chi_{k}\left(d\right)=i^{-k\left(1-d\right)}=i^{k\left(d-1\right)},$
which agrees with the quadratic residue symbol $\left(\frac{-1}{d}\right)^{k}$
if $d$ is odd. Note that $\chi_{k}$ is only a character on odd integers
but if $m=ad$ with $km\equiv\pm k\mod12$ then either $k$ is even
or $m$ is odd hence $\chi_{k}\left(d\right)$ is always real. 
\end{defn}

\subsection{Maass waveforms}

In this paper we fix an integer weight $k\not\equiv0\mod12$ (without
loss of generality $1\le k\le11$) and a multiplier system $v=v_{\eta}^{2k}$
given by $\xi=v\left(T\right)=e^{\frac{\pi ik}{6}}=e\left(\alpha\right),$
i.e. $\alpha=\frac{k}{12}$. Note that $\xi^{n}=1$ if $nk\equiv0\mod12$.
Set $D=\frac{k}{\left(12,k\right)}$. 

A Maass waveform on $\Gamma$ of weight $k$ and multiplier system
$v$ is a function $f\left(z\right)$ on $\H$ such that $f_{|\left[k,A\right]}=v\left(A\right)f\left(z\right)$
for all $A\in\overline{\Gamma}$, $\left(\Delta_{k}+\lambda\right)f=0$
with $\lambda=\frac{1}{4}+R^{2}\ge0$ and $\int_{\Gamma\backslash\H}\left|f\left(z\right)\right|^{2}d\mu\left(z\right)<\infty$
where $\mu\left(x+iy\right)=\frac{dxdy}{y^{2}}$ is the hyperbolic
area measure. Such $f$ is automatically cuspidal, i.e. $f\left(z\right)$
vanishes at $i\infty$ and to stress this some authors use the name
Maass cusp form instead of Maass waveform. For simplicity we only
consider $\lambda>\frac{1}{4}$, i.e. $R\in\mathbb{R}.$ The space
of all such Maass waveforms with these parameters is denoted by $\Mas{\Gamma,v,k,R}$
and the space of all Maass waveforms with multiplier $v$ and weight
$k$ and any spectral parameter is denoted by $\Mas{\Gamma,v,k}$.
To simplify certain Fourier expansions later we will use the notation
\begin{eqnarray*}
\mathcal{W}_{n,k,R}\left(z\right) & = & \frac{1}{\sqrt{\left|n+\alpha_{k}\right|}}W_{\frac{k}{2}\sgn\left(n+\alpha_{k}\right),iR}\left(4\pi\left|n+\alpha_{k}\right|y\right)e\left(\left(n+\alpha_{k}\right)x\right),\,\mbox{and}\\
\mathcal{W}_{n,k,R}^{*}\left(z\right) & = & \frac{1}{\sqrt{\left|n-\alpha_{k}\right|}}W_{\frac{k}{2}\sgn\left(n-\alpha_{k}\right),iR}\left(4\pi\left|n-\alpha_{k}\right|y\right)e\left(\left(n-\alpha_{k}\right)x\right),\end{eqnarray*}
for $z=x+iy\in\H$, where $W_{l,\mu}\left(y\right)$ denotes W-Whittaker
function in standard notation (cf. e.g.~\cite[vol.\ I, p.\ 264]{erdelyi:53})
and $\alpha_{k}=\frac{k}{12}$. It is known (cf. e.g.~\cite[ch.\ 9]{hejhal:lnm1001})
that any $f\in\Mas{\Gamma,v,k,R}$ has a Fourier expansion $f(z)=\sum_{n=-\infty}^{\infty}c\left(n\right)\mathcal{W}_{n,k,R}\left(z\right).$

\section{Main Result}

\begin{thm}
\label{thm:main_theorem}There exist a basis of $\Mas{\Gamma,v,k,R}$
consisting of Maass wave forms $f$ with Fourier expansions at infinity
\[
f(z)=\sum_{n=-\infty}^{\infty}c\left(n\right)\mathcal{W}_{n,k,R}\left(z\right)\]
where the coefficients $c\left(n\right)$ satisfy the following multiplicativity
relations if $c\left(0\right)\ne0$. For positive integers $m,n$
with $12m\equiv12n\equiv0\mod k$ set $m_{1}=\frac{12m}{k}$ and $n_{1}=\frac{12n}{k}$.
If $\left(m_{1}+1,D\right)=\left(n_{1}+1,D\right)=1$ then \begin{equation}
c(m)c(n)=c(0)\sum_{0<d|(m_{1}+1,n_{1}+1)}\chi_{k}\left(d\right)\, c\left(\frac{k}{12}\left(\frac{\left(m_{1}+1\right)\left(n_{1}+1\right)}{d^{2}}-1\right)\right),\label{eq:w1_coeff_basic_mult_rel}\end{equation}
whereas if $\left(m_{1}-1,D\right)=\left(n_{1}-1,D\right)=1$ then
\begin{equation}
c\left(-m\right)c\left(-n\right)=\Lambda_{k,R}\; c\left(0\right)\sum_{0<d|\left(m_{1}-1,n_{1}-1\right)}\chi_{k}\left(d\right)\, c\left(\frac{k}{12}\left(\frac{\left(m_{1}-1\right)\left(n_{1}-1\right)}{d^{2}}-1\right)\right),\label{eq:mult_rel_w1_neg_coef}\end{equation}
where \begin{equation}
\Lambda_{k,R}=\begin{cases}
\prod_{j=1}^{l}\left(j\left(j-1\right)+\frac{1}{4}+R^{2}\right)^{2}, & k=2l,\\
-R^{2}\prod_{j=1}^{l}\left(j^{2}+R^{2}\right)^{2}, & k=2l+1.\end{cases}\label{eq:Lambda_kR}\end{equation}

\end{thm}
\begin{proof}
We will prove this theorem in the following sections by constructing
appropriate self-adjoint Hecke (and Hecke-like) operators acting on
$\Mas{\Gamma,v,k,R}$ satisfying certain multiplicativity relations. 
\end{proof}
\begin{rem}
If $k|12$ then (\ref{eq:w1_coeff_basic_mult_rel}) and (\ref{eq:mult_rel_w1_neg_coef})
hold for all positive integers $m$ and $n$.
\end{rem}

\begin{rem}
The condition $c\left(0\right)\ne0$ is required to obtain this particular
form of relations (\ref{eq:w1_coeff_basic_mult_rel}) and (\ref{eq:mult_rel_w1_neg_coef}).
If $c\left(n_{0}\right)\ne0$ for some $n_{0}$ we get a similar set
of relations since the multiplicativity of the Hecke eigenvalues,
e.g.~(\ref{eq:mult_rel_lambda_mn}) and (\ref{eq:mult_rel_mu_mn})
will still hold while the explicit formulas (\ref{eq:lambda_n_explicit})
and (\ref{eq:mu_explicit}) have to be changed accordingly. If $\left(n_{0},D\right)=1$
the relations will even be as simple and if $k|12$ it is also easy
to show that if $f$ is a simultaneous Hecke eigenform then $c\left(0\right)\ne0$
unless $f\equiv0$. For the rest of the paper we assume for simplicity
that $c\left(0\right)\ne0$. 
\end{rem}
Explicit examples of multiplicativity relations are given in the last
section.

\section{Hecke Operators and Multiplier Systems}

The standard way to view Hecke operators today is via the action of
double cosets, cf.~e.g.~Shimura \cite{shimura} or Miyake \cite{miyake}.
It is possible to describe precise conditions for the existence of
Hecke operators defined by actions of double cosets also in the setting
of general multipliers and weights, cf.~\cite[§2.4.6]{stromberg:thesis}.
One problem that arises in the general setting is to make a consistent
definition of the operators $T(p,p)$ (in the terminology of \cite{shimura}
and \cite{miyake}). 

In part because of this difficulty we instead use a construction based
on the ideas of Wohlfahrt \cite{MR0106888} (also van Lint \cite{van-Lint:HeckeOperators:MR0090616})
which is more in the spirit of Hecke \cite{Hecke:37}. The key point
in this construction is to find the specific linear combination of
weight $k$ slash actions by elements of a complete set of $\Gamma$
inequivalent matrices of determinant $n$. The Hecke operators of
Wohlfahrt differs slightly from the standard construction for non-square
free integers. 

In this section we will encounter two different types of Hecke operators
$\T{m}$ which have different properties depending on whether $km\equiv\pm k\mod12$.
In the following section we construct a kind of {}``Hecke-like''
operator, $\Th{m}$, which consists of a Hecke operator $\T{m}$ composed
with another operator $\Theta$. These Hecke-like operators are crucial
in obtaining the new kind of multiplicativity relations (\ref{eq:mult_rel_w1_neg_coef}). 

\begin{prop}
\label{pro:def_general_Tm}Let $m$ be a positive integer satisfying
$km\equiv k\mod12$. We can then define a Hecke operator $\T{m}:\Mas{\Gamma,v,k,R}\rightarrow\Mas{\Gamma,v,k,R}$
by \begin{equation}
\T{m}f\left(z\right)=\frac{1}{\sqrt{n}}\sum_{ad=n,\, d>0}\chi_{k}\left(d\right)\sum_{b\mod d}\xi^{-bd}f\left(\frac{az+b}{d}\right)\label{eq:Hecke_n_general}\end{equation}
and if $f\left(z\right)$ has Fourier coefficients $c\left(n\right)$
then $\T{m}f\left(z\right)$ has Fourier coefficients $d\left(n\right)$
where \begin{eqnarray}
d\left(n\right) & = & \sum_{0<d|\left(m,n-k\frac{m-1}{12}\right)}\chi_{k}\left(d\right)c\left(\frac{nm}{d^{2}}+\frac{k\left(m-d^{2}\right)}{12d^{2}}\right).\label{eq:Formula_for_dn}\end{eqnarray}
If $\T{m}f=\lambda_{m,k}f$ and $c\left(0\right)\ne0$ then \begin{equation}
\lambda_{n,k}=\frac{1}{c\left(0\right)}\begin{cases}
c\left(k\frac{\left(n-1\right)}{12}\right), & \left(D,n\right)=1,\,\mbox{and}\\
c\left(k\frac{\left(n-1\right)}{12}\right)+\chi_{k}\left(D\right)c\left(k\frac{n-D^{2}}{12D^{2}}\right), & \left(D,n\right)>1.\end{cases}\label{eq:lambda_n_explicit}\end{equation}

\end{prop}
\begin{proof}
We follow the construction of Wohlfahrt \cite{MR0106888}. Let $Q=\left(\begin{smallmatrix}m & 0\\
0 & 1\end{smallmatrix}\right)$, $\beta_{a,b,d}=\left(\begin{smallmatrix}a & b\\
0 & d\end{smallmatrix}\right)$ and define the sets \begin{eqnarray*}
\mathcal{R}_{\, m} & = & \left\{ \beta_{a,b,d};\, ad=m,\, b\mod d,\,\left(a,b,d\right)=1\right\} ,\,\mbox{and}\\
\mathcal{R}_{\, m}^{*} & = & \left\{ \beta_{a,b,d};\, ad=m,\, b\mod d\right\} .\end{eqnarray*}
It is well-known (cf.~e.g.~\cite[p.\ 244]{MR0106888}, \cite[p.\ 142]{miyake}
or \cite[p.\ 167]{koblitz-introduction-to-elliptic}) that both $\mathcal{R}_{\, m}$and
$\mathcal{R}_{\, m}^{*}$ are complete sets of left $\Gamma$-inequivalent
matrices in $\GLZ$ with determinant $n$. The Wohlfahrt operator
$T\left\langle Q\right\rangle $ given by \begin{equation}
T\left\langle Q\right\rangle f=\frac{1}{\sqrt{m}}\sum_{\beta\in\mathcal{R}_{\, m}}\overline{v_{Q}\left(\beta\right)}\, f_{|[k,\beta]}\left(z\right),\label{eq:Wholfahrt_TQ}\end{equation}
where $v_{Q}\left(\beta\right)=v_{1}\left(A\right)v_{2}\left(B\right)$
for $\beta=AQB\in\mathcal{R}_{\, m}$ is shown to map the space $\Mas{\Gamma,v_{1},k,R}$
into $\Mas{\Gamma,v_{2},k,R}$. In our case, with $Q$ as above $v_{1}=v$
and $v_{2}=v^{m}=v$ since $mk\equiv k\mod12$. To compute the numbers
$v_{Q}\left(\beta\right)$ explicitly we use the explicit formula
(\ref{eq:v_eta_explicit}) together with $v\left(-A\right)=e^{-\pi ik}v\left(A\right)$
valid for all $A\in\overline{\Gamma}$ and the decomposition $\left(\begin{smallmatrix}a & b\\
0 & d\end{smallmatrix}\right)=\left(\begin{smallmatrix}0 & 1\\
-1 & du\end{smallmatrix}\right)Q\left(\begin{smallmatrix}u & -v\\
a & b\end{smallmatrix}\right),$ $ub+va=1$ . We now see that \begin{eqnarray*}
v_{Q}\left(\beta_{a,b,d}\right) & = & v\left(\begin{smallmatrix}0 & 1\\
-1 & du\end{smallmatrix}\right)v\left(\begin{smallmatrix}u & -v\\
a & b\end{smallmatrix}\right)\\
 & = & \xi^{bd}\xi^{3\left(1-d\right)}=\xi^{bd}i^{k\left(1-d\right)}=\xi^{bd}\chi_{k}\left(d\right),\end{eqnarray*}
where $\chi_{k}\left(d\right)=\left(\frac{-1}{d}\right)^{k}$ for
odd $d$ and is in our cases always real since either $k$ is even
or $d$ is odd. From (\ref{eq:Hecke_n_general}) we see that $\T{m}$
is defined by a similar sum but instead over $\mathcal{R}_{\, n}^{*}$.
Clearly \begin{eqnarray*}
\mathcal{R}_{\, m}^{*}\smallsetminus\mathcal{R}_{\, m} & = & \left\{ \beta_{la',lb',ld'}=\left(\begin{smallmatrix}a'l & b'l\\
0 & d'l\end{smallmatrix}\right);\, l^{2}|m,\, b'\mod d'\right\} \\
 & = & \cup_{l^{2}|m}\mathcal{R}_{\,\frac{m}{l^{2}}}\end{eqnarray*}
and thus \[
\T{m}f\left(z\right)=T\left\langle Q\right\rangle +\sum_{l^{2}|m}\frac{1}{l}\chi_{k}\left(l\right)\T{\frac{m}{l^{2}}}.\]
Using induction on the number of prime factors of $m$ together with
the properties of $T\left\langle Q\right\rangle $ we see that $\T{n}$
is indeed a map from $\Mas{\Gamma,v,k,R}$ to $\Mas{\Gamma,v,k,R}.$
The Fourier coefficients are easily obtained by direct computation:
\begin{eqnarray*}
\T{m}f & = & \frac{1}{\sqrt{m}}\sum_{ad=m,\, d>0}\chi_{k}\left(d\right)\sum_{b\mod d}\xi^{-bd}\sum_{n}c\left(n\right)\mathcal{W}_{n,k,R}\left(\frac{az+b}{d}\right)\\
 & = & \frac{1}{\sqrt{m}}\sum_{ad=m,\, d>0}\chi_{k}\left(d\right)\sum_{n}c\left(n\right)\mathcal{W}_{n,k,R}\left(\frac{az}{d}\right)\sum_{b=0}^{d-1}e\left(\frac{b}{12}\left[-dk+\left(n+\frac{k}{12}\right)\frac{1}{d}\right]\right)\end{eqnarray*}
and using the Gauss sum formula\[
\sum_{b\mod d}e\left(\frac{b}{12d}\left[12n+k-bd^{2}k\right]\right)=\begin{cases}
d, & \mbox{if }n+\frac{k}{12}=d\left(l+\frac{bdk}{12}\right),\, l\in\mathbb{Z},\\
0, & \mbox{else}\end{cases}\]
it is easy to see that\begin{eqnarray*}
\T{m}f\left(z\right) & = & \frac{1}{\sqrt{m}}\sum_{ad=m}\chi_{k}\left(d\right)\sqrt{d}\sum_{l=-\infty}^{\infty}c\left(dl+\frac{k\left(d^{2}-1\right)}{12}\right)\mathcal{W}_{al+\frac{k\left(m-1\right)}{12},k,R}\left(z\right)\\
 & = & \sum_{ad=m}\chi_{k}\left(d\right)\sum_{n\equiv\frac{k\left(m-1\right)}{12}\mod a}c\left(\frac{mn}{a^{2}}+\frac{k\left(m-a^{2}\right)}{12a^{2}}\right)\mathcal{W}_{n,k,R}\left(z\right)\\
 & = & \sum_{n=-\infty}^{\infty}\left\{ \sum_{0<a|\left(m,n-\frac{k\left(m-1\right)}{12}\right)}\chi_{k}\left(a\right)c\left(\frac{mn}{a^{2}}+\frac{k\left(m-a^{2}\right)}{12a^{2}}\right)\right\} \mathcal{W}_{n,k,R}\left(z\right).\end{eqnarray*}
The formula for the eigenvalues $\lambda_{m,k}$ are obtained by comparing
the coefficients $c\left(0\right)$ and $d\left(0\right)$ in (\ref{eq:Formula_for_dn})
and observing that $\left(m,\frac{k\left(m-1\right)}{12}\right)=\left(m,\frac{k}{\left(12,k\right)}\right)=\left(m,D\right)$.
The multiplicativity relation follows from Lemma \ref{lem:Heckeop_mult_rel}
below.
\end{proof}
\begin{rem}
In our case, if $km\equiv\pm k\mod12$ and $\left(m,D\right)>1$ it
can be shown that\\
$m\equiv\pm D^{2}\mod\frac{12D}{\left(12,k\right)}$. 
\end{rem}

\begin{rem}
By changing the action $f_{|[k,\beta]}$ to the holomorphic weight
$k$ action, $F_{|\left\{ k,\beta\right\} },$ the Hecke operators
$\T{m}$ for $km\equiv k\mod12$ can also be defined for holomorphic
modular forms of weight $k$ and multiplier system $v$.
\end{rem}
\begin{prop}
\label{pro:Second_type_Heckeoperator}If $m$ is a positive integer
and $km\equiv-k\mod12$ we define two different Hecke operators \begin{eqnarray*}
\T{m}:\Mas{\Gamma,v,k,R} & \rightarrow & \Mas{\Gamma,\overline{v},k,R},\,\mbox{and}\\
\T[\overline{v}]{m}:\Mas{\Gamma,\overline{v},k,R} & \rightarrow & \Mas{\Gamma,v,k,R}.\end{eqnarray*}
If $f\in\Mas{\Gamma,v,k,R}$ and $g\in\Mas{\Gamma,\overline{v},k,R}$
have Fourier expansions $f\left(z\right)=\sum_{n}c\left(n\right)\mathcal{W}_{n,k,R}\left(z\right)$
and $g\left(z\right)=\sum_{n}a\left(n\right)\mathcal{W}_{n,k,R}^{*}\left(z\right)$
respectively then \begin{eqnarray*}
\T{m}f\left(z\right) & = & \sum_{n=-\infty}^{\infty}d\left(n\right)\mathcal{W}_{n,k,R}^{*}\left(z\right),\,\mbox{and}\\
\T[\overline{v}]{m}g\left(z\right) & = & \sum_{n=-\infty}^{\infty}b\left(n\right)\mathcal{W}_{n,k,R}\left(z\right),\end{eqnarray*}
where \begin{eqnarray*}
d\left(n\right) & = & \sum_{0<a|\left(m,n-\frac{k\left(m+1\right)}{12}\right)}\chi_{k}\left(a\right)c\left(\frac{mn}{a^{2}}-\frac{k\left(m+a^{2}\right)}{12a^{2}}\right),\mbox{ and}\\
b\left(n\right) & = & \sum_{0<a|\left(m,n+\frac{k\left(m+1\right)}{12}\right)}\chi_{k}\left(a\right)c\left(\frac{mn}{a^{2}}+\frac{k\left(m+a^{2}\right)}{12a^{2}}\right).\end{eqnarray*}

\end{prop}
\begin{proof}
The construction of these two operators are exactly as above with
the exception that now in (\ref{eq:Wholfahrt_TQ}) we have $v_{Q}\left(AQB\right)=v_{1}\left(A\right)v_{2}\left(B\right)$
with $v_{1}=v$ and $v_{2}=v^{m}=\overline{v}$ for $\T{m}$ and $v_{1}=\overline{v}$
and $v_{2}=\overline{v}^{m}=v$ for $\T[\overline{v}]{m}$. $\T{m}$
is then given by the same formula as (\ref{eq:Hecke_n_general}),
while $\T[\overline{v}]{m}$ is given by the same formula with $\xi$
interchanged with $\overline{\xi}=\xi^{-1}$. 
\end{proof}
\begin{rem}
Note that $\T{1}$ can always be defined as the identity operator.
But in the case $k=6$, i.e. $k\equiv-k\mod12$ it can also be defined
as in Proposition \ref{pro:Second_type_Heckeoperator}, still acting
as the identity but shifting the summation in the Fourier series ($n\rightarrow n-1$). 
\end{rem}

\begin{rem}
Observe that the difference between the sets $\mathcal{R}_{m}$ and
the sets $\mathcal{R}_{m}^{*}$ we are actually summing over in $\T{m}$consist
of certain {}``trivial'' elements of lower order. In the case $m=p^{2}$
the trivial element is $\beta_{0}=\left(\begin{smallmatrix}p & 0\\
0 & p\end{smallmatrix}\right)$ which acts as the identity in $\PSLR$. When the multiplier is trivial
or a Dirichlet character $\chi$ this element acts as multiplication
by the scalar $\chi\left(p\right)$ (in the notation of \cite{miyake,shimura}
this operator is denoted by $T$$\left(p,p\right)$).

The question of including or excluding these lower order elements
is precisely what distinguishes the classical Hecke operators for
integral weights and characters or half integral weights and the theta
multiplier system on one hand and the general Wohlfahrt operators
on the other hand. Cf.~e.g.~\cite{shimura:73:half_integral} and
also \cite[p.\ 206]{koblitz-introduction-to-elliptic} where the difference
is mentioned.

In practice, the exclusion of the trivial elements results in more
intricate multiplicativity relations. Cf. e.g. \cite[\S 9]{MR0106888}.
In relation with this discussion it should also be mentioned that
the {}``Hecke like'' operators of \cite{MR2155506,MR2180234} (acting
on on period functions) correspond precisely to the exclusion of $\beta_{0}$
in the definition of the Hecke operator for $p^{2}$ (compare the
set $X_{n}^{*}$ in \cite[p.\ 143]{MR2180234} with the set $\mathcal{R}_{n}$). 

We are now in a position to start obtaining the multiplicativity relations.
The following Proposition is crucial. 
\end{rem}
\begin{prop}
\label{lem:Heckeop_mult_rel}For any \emph{positive} integers $n,m$,
$kn,km\equiv\pm k\mod12$ we have \begin{eqnarray*}
\T{m}\T{n} & = & \sum_{0<d|\left(m,n\right)}\chi_{k}\left(d\right)\T{\frac{mn}{d^{2}}},\, kn\equiv k\mod12,\\
\T[\overline{v}]{m}\T{n} & = & \sum_{0<d|\left(m,n\right)}\chi_{k}\left(d\right)\T{\frac{mn}{d^{2}}},\, kn\equiv-k\mod12.\end{eqnarray*}

\end{prop}
The following three lemmas are easy to verify directly by using the
definitions and comparing both sides of the equalities. 

\begin{lem}
\label{lem:Heckeop_comm_rel}For any primes $p_{1}\ne p_{2},q_{1}\ne q_{2}\ge3$,
$kp_{j}\equiv-kq_{i}\equiv k\mod12$ the following commutation relations
hold \begin{eqnarray*}
\T{p_{1}}\T{p_{2}} & = & \T{p_{2}}\T{p_{1}}=\T{p_{1}p_{2}},\\
\T[\overline{v}]{q_{1}}\T{q_{2}} & = & \T[\overline{v}]{q_{2}}\T{q_{1}}=\T{q_{1}q_{2}}.\end{eqnarray*}
Furthermore, for primes $p,q\ge3$, $kp\equiv-kq\equiv k\mod12$ we
have \[
\T{q}\T{p}=\T[\overline{v}]{p}\T{q}=\T{pq}.\]

\end{lem}

\begin{lem}
\label{lem:Hecek_op_mult_one_prime}For any pair of primes $p,q\ge3$
with $kp\equiv-kq\equiv k\mod12$ we have \begin{eqnarray*}
\T{p^{r}}\T{p} & = & \T{p^{r+1}}+\chi_{k}\left(p\right)\T{p^{r-1}},\, r\ge1,\\
\T[\overline{v}]{q^{r}}\T{q} & = & \T{q^{r+1}}+\chi_{k}\left(q\right)\T{q^{r-1}},\, r\ge1.\end{eqnarray*}

\end{lem}
By induction one then deduce the following Lemma. 

\begin{lem}
\label{lem:Heckeop_mult_rel_more}For any primes $p,q\ge2$ and pair
of integers $s,r\ge1$ we have \begin{eqnarray*}
\T{p^{r}}\T{p^{s}} & = & \sum_{d|\left(p^{r},p^{s}\right)}\chi_{k}\left(d\right)\T{\frac{p^{s+r}}{d^{2}}},\,\mbox{if }kp^{r}\equiv kp^{s}\equiv k\mod12,\\
\T[\overline{v}]{q^{r}}\T{q^{s}} & = & \sum_{d|\left(q^{r},q^{s}\right)}\chi_{k}\left(d\right)\T{\frac{q^{s+r}}{d^{2}}},\,\mbox{if }kq^{r}\equiv kq^{s}\equiv-k\mod12,\, s\,\mbox{odd,}\\
\T{q^{r}}\T{q^{s}} & = & \sum_{d|\left(q^{r},q^{s}\right)}\chi_{k}\left(d\right)\T{\frac{q^{s+r}}{d^{2}}},\,\mbox{if }kq^{r}\equiv kq^{s}\equiv-k\mod12,\, s\,\mbox{even.}\end{eqnarray*}

\end{lem}
\begin{proof}[Proof of Proposition \ref{lem:Heckeop_mult_rel}]

The proposition now follows by induction from Lemmas \ref{lem:Heckeop_comm_rel}-\ref{lem:Heckeop_mult_rel_more}.\end{proof}

\begin{lem}
\label{lem:mult_rel_prime_Tm_Tn}If $\left(m,n\right)=1$ then \begin{eqnarray*}
\T{mn} & = & \T{m}\T{n},\, km\equiv kn\equiv k\mod12,\\
\T{mn} & = & \T[\overline{v}]{m}\T{n},\, km\equiv kn\equiv-k\mod12,\\
\T{mn} & = & \T[\overline{v}]{m}\T{n}=\T{n}\T{m},\, km\equiv-kn\equiv k\mod12.\end{eqnarray*}

\end{lem}
\begin{proof}
Together with the definitions, use the relation $\beta_{a_{1},b_{1},d_{1}}\beta_{a_{2},b_{2},d_{2}}=\beta_{a,b,d}$
with $a=a_{1}a_{2},$ $d=d_{1}d_{2}$ and $b=a_{1}b_{2}+d_{2}b_{1}$
and observe that $b$ runs through a complete set of residues modulo
$d$ as $b_{1}$ and $b_{2}$ runs through residues of $d_{1}$ and
$d_{2}$ respectively. 
\end{proof}

It is clear that as usual the Hecke operators $\T{p}$ with $p$ prime
generate the complete algebra and that the Hecke eigenvalues satisfy
the desired multiplicativity relations. 

\begin{prop}
\label{pro:mult_rel_Tm}If $f(z)\in\Mas{\Gamma,v,k,R}$ is an eigenfunction
of all $\T{m},$ $km\equiv k\mod12$ with eigenvalues $\lambda_{m}$
then the following multiplicativity relation holds for all positive
integers $m,n$, $km\equiv kn\equiv k\mod12$ \begin{equation}
\lambda_{m}\lambda_{n}=\sum_{d|\left(m,n\right)}\chi_{k}\left(d\right)\lambda_{\frac{mn}{d^{2}}}.\label{eq:mult_rel_lambda_mn}\end{equation}
Furthermore, if $f\left(z\right)=\sum_{n}c\left(n\right)\mathcal{W}_{n,k,R}\left(z\right)$
and if $\left(m,D\right)=\left(n,D\right)=1$ then \[
c\left(k\frac{m-1}{12}\right)c\left(k\frac{n-1}{12}\right)=c\left(0\right)\sum_{d|\left(m,n\right)}\chi_{k}\left(d\right)c\left(k\frac{mn-d^{2}}{12d^{2}}\right).\]

\end{prop}
From the last lemma we obtain the relation (\ref{eq:w1_coeff_basic_mult_rel})
of Theorem \ref{thm:main_theorem}. In Section \ref{sub:The-auxillary-operator}
we will obtain multiplicativity relations for the negative coefficients
by introducing an auxiliary operator $\Theta$ which composed with
$\T{q}$ preserves the multiplier system. 

\begin{defn}
For any pair of functions $f,g\in\Mas{\Gamma,v,k}$ or $\Mas{\Gamma,\overline{v},k}$
we define the \emph{Petersson inner product} (cf.~\cite{MR0060542})
by \[
\left(f,g\right)=\int_{\Gamma\backslash\H}f\left(z\right)\overline{g\left(z\right)}y^{k-2}dxdy.\]

\end{defn}
\begin{lem}
\label{lem:relation_Tn*andTn}Let $n$ and $m$ be positive integers
with $kn\equiv-km\equiv k\mod12$ and let $\T{n},$ $\T{m}$ and $\T[\overline{v}]{m}$
be the Hecke operators defined in Propositions \ref{pro:def_general_Tm}
and \ref{pro:Second_type_Heckeoperator}. Let $\Ta{n}$, $\Ta{m}$
and $\Ta[\overline{v}]{m}$ be the respective adjoint operators with
respect to the Petersson inner product. Then \begin{eqnarray*}
\Ta{n} & = & \T{n},\\
\Ta{m} & = & \left(-1\right)^{k}\T[\overline{v}]{m},\quad\mbox{and}\\
\Ta[\overline{v}]{m} & = & \left(-1\right)^{k}\T{m}.\end{eqnarray*}

\end{lem}
\begin{proof}
We can use the expression for the adjoint of the Wohlfahrt operator
given in \cite[§4 and §6]{MR0106888}\[
T\left\langle Q\right\rangle ^{*}=\sum_{\beta\in\mathcal{R}_{m}}\overline{v_{Q}^{*}\left(\beta\right)}\, f\left(\beta z\right)\]
where $\overline{v_{Q}^{*}\left(\beta\right)}=v_{Q}\left(\beta^{*}\right)$,
$\beta=\beta_{a,b,d}=\left(\begin{smallmatrix}a & b\\
0 & d\end{smallmatrix}\right)$ and $\beta^{*}=\left(\begin{smallmatrix}d & -b\\
0 & a\end{smallmatrix}\right)$. The identification of the adjoints above with the right hand sides
are thus easily done by evaluating $v_{Q}\left(\beta^{*}\right)$
(as in the proof of Proposition \ref{pro:def_general_Tm}) and comparing
with the corresponding factors in $\T{n},$$\T[\overline{v}]{m}$
and $\T{m}$ respectively. 
\end{proof}
\begin{cor}
If $f\in\Mas{\Gamma,v,k}$ is an eigenfunction of all $\T{n}$ with
$kn\equiv k\mod12$ and has Fourier coefficients $c\left(l\right)$
then the quotients $\frac{1}{c\left(0\right)}c\left(l\right)$ are
real for all $l$ such that $k|12l$. 
\end{cor}

\subsection{\label{sub:The-auxillary-operator}The auxiliary operator $\Theta$}

The operator $\Theta$ is essentially the operator defined by Maass
in \cite[p.\ 181]{maass:modular_functions}. $\Theta$ interchanges
the multiplier system with its conjugate and while Maass used $\Theta$
in connection with a real multiplier system we will compose it with
operators $\T{m}$ for $km\equiv-k\mod12$ to obtain Hecke-like operators
preserving the multiplier system.

Let \begin{eqnarray*}
E_{k}^{+}=(z-\overline{z})\frac{\partial}{\partial z}+\frac{k}{2}, & \,\textrm{and } & E_{k}^{-}=-(z-\overline{z})\frac{\partial}{\partial\overline{z}}+\frac{k}{2}\end{eqnarray*}
be the Maass raising and lowering operators. Cf.~e.g.~\cite[p.\ 369 and pp.\ 381-382]{hejhal:lnm1001}
or \cite[p.\ 188]{maass:modular_functions}. It is known that for
any multiplier system $v$ of weight $k$ the operator $E_{k}^{\pm}$
is one-to-one and onto from $\Mas{\Gamma,v,k,R}$ to $\Mas{\Gamma,v,k\pm2,R}$
if $\lambda=\frac{1}{4}+R^{2}\ne\mp\frac{k}{2}\left(1\pm\frac{k}{2}\right)$.

We now define $\mathcal{E}_{k}^{-}=E_{2-k}^{-}\circ\cdots\circ E_{k-2}^{-}\circ E_{k}^{-}$
and $\mathcal{E}_{-k}^{+}=E_{k-2}^{+}\circ\cdots\circ E_{-k+2}^{+}\circ E_{-k}^{+}$.
For $\lambda>\frac{1}{4}$ it is clear that $\mathcal{E}_{k}^{-}$
is a one-to-one and onto map from $\Mas{\Gamma,v,k,R}$ to $\Mas{\Gamma,v,-k,R}$
and that $\mathcal{E}_{-k}^{+}$ is one-to-one and onto in the other
direction. 

Let $J$ denote the reflection in the imaginary axis, i.e. $Jf(z)=f\left(-\overline{z}\right)$.
Clearly $J$ is an involution and by using Lemma \ref{lem:mult_system_by_xi}
and it is also easy to see that $J$ maps $\Mas{\Gamma,v,k,R}$ to
$\Mas{\Gamma,\overline{v},-k,R}$ (cf. also \cite[§2.4.1]{stromberg:thesis}).

\begin{defn}
We can now define the operator $\Theta$ as the composition of $J$
and $\mathcal{E}_{k}^{-}$: \[
\Theta=J\mathcal{E}_{k}^{-}.\]
This operator is clearly one-to-one and onto from $\Mas{\Gamma,v,k,R}$
to $\Mas{\Gamma,\overline{v},k,R}$ as well as from $\Mas{\Gamma,\overline{v},k,R}$
to $\Mas{\Gamma,v,k,R}$. 
\end{defn}
The following properties of $\Theta$ are now easily verified.

\begin{lem}
For any $f\in\Mas{\Gamma,v,k,\frac{1}{4}+R^{2}}$ we have \[
\Theta^{2}f=\Lambda_{k,R}f,\]
where $\Lambda_{k,R}$ is given by (\ref{eq:Lambda_kR}).
\end{lem}
\begin{proof}
The action of the Maass operators on a Fourier series is given by
\[
E_{k}^{\pm}\left[W_{\frac{k}{2}\epsilon,iR}\left(y\right)e\left(x\right)\right]=\pm W_{\frac{k\pm2}{2}\epsilon}\left(y\right)e\left(x\right)\times\begin{cases}
-1, & \epsilon=\pm1,\\
\left(\frac{k\left(k\pm2\right)}{4}+\frac{1}{4}+R^{2}\right), & \epsilon=\mp1.\end{cases}\]
Using this formula it is clear that if $f\left(z\right)=\sum_{n}c\left(n\right)\mathcal{W}_{n,k}\left(z\right)$
then \[
\Theta f=J\mathcal{E}_{k}^{-}f(z)=\sum_{n=-\infty}^{\infty}\delta_{n}c\left(-n\right)\mathcal{W}_{n,k,R}^{*}\left(z\right),\]
where $\delta_{n}=1$ if $n\ge1$ and $\delta_{n}=\Lambda_{k,R}$
if $n\le0$. That \begin{equation}
\Lambda_{k,R}=\left(-1\right)^{k}\prod_{j=0}^{k-1}\left(\frac{\left(k-2j\right)\left(k-2-2j\right)+1}{4}+R^{2}\right)\label{eq:Lambda_kr_first_ver}\end{equation}
is also given by the formula of theorem is easy to verify for all
integers $k,$ $1\le k\le11.$ It is now easy to see that $\Theta^{2}f=\Lambda_{k,R}f$
.
\end{proof}
\begin{lem}
The adjoint of $\Theta$ with respect to the Petersson inner product
is given by \[
\Theta^{*}=\left(-1\right)^{k}\Theta.\]

\end{lem}
The formula of Theorem \ref{thm:main_theorem} is then easy to verify
for all $k$. The lowering and raising operators are adjoints with
respect to the Petersson inner product in the sense that $E_{k}^{-*}=E_{k-2}^{+}$
(cf. e.g. \cite[pp.\ 135-136]{bump} and note that $E_{k}^{+}=R_{k}$
while $E_{k}^{-}=-L_{k}$). It is clear that the adjoint of $\Theta$
is given by \[
\Theta^{*}=\mathcal{E}_{-k}^{+}J\]
where $\mathcal{E}_{-k}^{+}$ is the adjoint of $\mathcal{E}_{k}^{-}$.
By direct verification it is now easy to see that \[
\Theta^{*}=\left(-1\right)^{k}\Theta.\]

\begin{defn}
For any positive integer $m$ with $km\equiv-k\mod12$ we define \[
\Th{m}=\T[\overline{v}]{m}\Theta.\]

\end{defn}
\begin{lem}
\label{lem:Theta_commutes_with_Tm}The operator $\Theta$ commutes
with the Hecke operators, i.e. for $km\equiv-k\mod12$ \[
\Th{m}=\T[\overline{v}]{m}\Theta=\Theta\T{m}.\]

\end{lem}
\begin{proof}
By direct verification. 
\end{proof}
\begin{lem}
\label{lem:Th_self_adjoint}$\Th{m}$ is self-adjoint for any positive
integer $m$ with $km\equiv-k\mod12$. 
\end{lem}
\begin{proof}
This follows clearly from the fact that $\Theta^{*}=\left(-1\right)^{k}$$\Theta$
and $\Ta{m}=\left(-1\right)^{k}\T[\overline{v}]{m}$.
\end{proof}
\begin{prop}
\label{pro:Th{m}_action}Let $m$ be a positive integer with $km\equiv-k\mod12$
we can then define the Hecke-like operator $\Th{m}=\T[\overline{v}]{m}\Theta$
on $\Mas{\Gamma,v,k,R}$. Let $f(z)\in\Mas{\Gamma,v,k,R}$ be given
by the Fourier expansion $f\left(z\right)=\sum_{n}c\left(n\right)\mathcal{W}_{n,k,R}\left(z\right)$
then\[
\Th{m}f(z)=\sum_{n=-\infty}^{\infty}d\left(n\right)\mathcal{W}_{n,k,R}\left(z\right),\]
where\begin{eqnarray*}
d\left(n\right) & = & \sum_{0<a|\left(m,n+\frac{k\left(m+1\right)}{12}\right)}\chi_{k}\left(a\right)\delta_{a}c\left(-\frac{mn}{a^{2}}-\frac{k\left(m+a^{2}\right)}{12a^{2}}\right),\\
 &  & \delta_{a}=\begin{cases}
1, & \frac{mn}{a^{2}}+\frac{k\left(m+a^{2}\right)}{12a^{2}}\ge1,\\
\Lambda_{k,R}, & \mbox{else.}\end{cases}\end{eqnarray*}
In particular, if $\Th{m}f=\mu_{m}f$ and $c\left(0\right)\ne0$ then
\begin{equation}
\mu_{m}=\frac{1}{c\left(0\right)}\begin{cases}
c\left(-k\frac{1+m}{12}\right), & \left(m,D\right)=1,\mbox{ and}\\
c\left(-k\frac{1+m}{12}\right)+\chi_{k}\left(D\right)c\left(-k\frac{m+D^{2}}{12D^{2}}\right), & \left(m,D\right)>1.\end{cases}\label{eq:mu_explicit}\end{equation}

\end{prop}
\begin{proof}
The formula for $d\left(n\right)$ is obtained by direct computation.
For the eigenvalue, $\mu_{m},$ a similar argument as in the proof
of Proposition \ref{pro:def_general_Tm} shows that $\left(m,\frac{k\left(m+1\right)}{12}\right)=\left(m,D\right)$
and $D$ is necessarily prime since $k\in[1,11]\cap\mathbb{Z}$. 
\end{proof}
\begin{prop}
\label{pro:mult_rel_Th}For any positive integers $m_{1}$and $m_{2}$
such that $km_{1}\equiv km_{2}\equiv-k\mod12$ \[
\Th{m_{1}}\Th{m_{2}}=\Lambda_{k,R}\sum_{d|\left(m_{1},m_{2}\right)}\chi_{k}\left(d\right)\T{\frac{m_{1}m_{2}}{d^{2}}}.\]
We deduce that if $f$ is an eigenform of $\Th{q}$ and all $\T{p}$
for all primes with $kp\equiv-kq\equiv k\mod12$ it is an eigenform
of all $\T{n}$ and all $\Th{m}$ with $kn\equiv-km\equiv k\mod12$.
If we denote the eigenvalues by $\lambda_{n}$ and $\mu_{m}$ respectively
we have \begin{equation}
\mu_{m}\mu_{n}=\Lambda_{k,R}\sum_{d|\left(m,n\right)}\chi_{k}\left(d\right)\lambda_{\frac{mn}{d^{2}}}.\label{eq:mult_rel_mu_mn}\end{equation}
If $f\left(z\right)=\sum_{n}c\left(n\right)\mathcal{W}_{n,k,R}\left(z\right)$
then, for all $km\equiv kn\equiv-k\mod12$ with $\left(m,D\right)=\left(n,D\right)=1$
\[
c\left(-k\frac{m+1}{12}\right)c\left(-k\frac{n+1}{12}\right)=\Lambda_{k,R}c\left(0\right)\sum_{d|\left(m,n\right)}\chi_{k}\left(d\right)c\left(k\frac{mn-d^{2}}{12d^{2}}\right).\]

\end{prop}
\begin{proof}
The first part is easily established by direct verification. The multiplicativity
relation follows from Lemma \ref{lem:Heckeop_mult_rel} and \ref{lem:Theta_commutes_with_Tm}
and (\ref{eq:mult_rel_w1_neg_coef}) is now immediate by interchanging
$m$ with $k\frac{m+1}{12}$.
\end{proof}
\begin{lem}
The family $\mathcal{T}=\left\{ \T{n}\right\} _{kn\equiv k\mod12}\cup\left\{ \Th{m}\right\} _{km\equiv-k\mod12}$
consists of self-adjoint operators on $\Mas{\Gamma,v,k}$. Furthermore
the operators in $\mathcal{T}$ also commute amongst themselves and
with the weight $k$ Laplacian $\Delta_{k}$. 
\end{lem}
\begin{proof}
The self-adjointness follows from Lemmas \ref{lem:relation_Tn*andTn}
and \ref{lem:Th_self_adjoint}. That the operators commute with each
other other follows from Proposition \ref{lem:Heckeop_mult_rel} and
Lemma \ref{lem:Theta_commutes_with_Tm}. $\Delta_{k}$ commutes with
the Hecke operators due to the well-known fact that it commutes with
the weight $k$ slash-action. To verify this fact it is easiest to
write $\Delta_{k}$ and the slash action in terms of $z$ and $\overline{z}$
using $\frac{\partial}{\partial z}=\frac{1}{2}\left[\frac{\partial}{\partial x}-i\frac{\partial}{\partial y}\right]$,
$\frac{\partial}{\partial\overline{z}}=\frac{1}{2}\left[\frac{\partial}{\partial x}+i\frac{\partial}{\partial y}\right]$
and $e^{i\Arg\left(cz+d\right)}=\left(\frac{cz+d}{c\overline{z}+d}\right)^{\frac{1}{2}}$. 
\end{proof}
A standard theorem from linear algebra about simultaneously diagonalization
of self-adjoint linear operators (although normality is enough) now
proves the following corollary. Together with Propositions \ref{pro:mult_rel_Tm}
and \ref{pro:mult_rel_Th} this concludes the proof of Theorem \ref{thm:main_theorem}.

\begin{cor}
There exist a basis of $\Mas{\Gamma,v,k,R}$ consisting of eigenfunctions
of all $\T{n}$ and $\Th{m}$ for $kn\equiv-km\equiv k\mod12$. 
\end{cor}

\section{Numerical Examples}

In this section we will write down some examples of multiplicativity
relations for various weights $k$, chosen so as to illustrate all
relevant properties. Using algorithms for computing Maass waveforms
with non-trivial multiplier systems as detailed in e.g. chapter two
of \cite{stromberg:thesis} we also demonstrate some numerical evidence
for these multiplicativity relations.

Any $f\in\Mas{\Gamma,v,k,R}$ have a Fourier expansion $f\left(z\right)=\sum_{n}c\left(n\right)\mathcal{W}_{n,k}\left(z\right)$
and we use either the normalization $c\left(1\right)=1$ or $c\left(0\right)=1$.
If nothing else is stated $c\left(1\right)=1$ but in certain cases
(e.g. $k=6$) the normalization $c\left(0\right)=1$ results in much
more stable numerics. We will also use the standard notation that
$c\left(x\right)=0$ if $x$ is not an integer.

Observe that this type of arithmetical normalization implies that
the relative sizes of coefficients with positive and negative indices
respectively differ with several orders of magnitude. This can most
easily be seen in the relations like (\ref{eq:mult_rel_w1_neg_coef})
relating on one side negative coefficients and on the other side positive
coefficients. In such relations the factor $\sqrt{\left|\Lambda_{k,R}\right|}$
is clearly distinguishable as a good measure of the difference in
magnitude. 

\begin{example}
Let $k=2,$ then $\left(k,12\right)=2$ and $D=1$. For $m\equiv-n\equiv1\mod6$
we have \[
\lambda_{m}=\frac{1}{c\left(0\right)}c\left(\frac{m-1}{6}\right),\,\mu_{n}=\frac{1}{c\left(0\right)}c\left(-\frac{m+1}{6}\right),\,\mbox{and}\,{\,\Lambda}_{2,R}=\left(\frac{1}{4}+R^{2}\right)^{2}.\]
For $m=7$ and $n=13$ we have $\lambda_{7,2}\lambda_{13,2}=\lambda_{91,2}$
or equivalently \[
c\left(1\right)c\left(2\right)=c\left(0\right)c\left(15\right).\]
For $m=n=7$ we have $\lambda_{7,2}^{2}=\lambda_{49,2}+\chi_{2}\left(7\right)\lambda_{1}$
or equivalently \[
c\left(1\right)^{2}=c\left(0\right)\left[c\left(8\right)+c\left(0\right)\right].\]
For $m=5$ and $n\equiv-1\mod6$ and $\left(5,n\right)=1$ we have
\[
\mu_{5,2}\mu_{n,2}=\Lambda_{k,R}\lambda_{5n,2},\]
or equivalently, using a change of variables to $l=\frac{n+1}{6}$:
\begin{equation}
c\left(1\right)c\left(-l\right)=c\left(0\right)\Lambda_{k,R}\left[c\left(5l-1\right)+c\left(\frac{l-1}{5}\right)\right],\label{eq:rel_w2_cl}\end{equation}
Table \ref{tab:N1w2R295} contains numerical values of Fourier coefficients
for $f\in\Mas{\Gamma,v,2.95645894117486}$. The third column of this
table consists of numerical values of (\ref{eq:rel_w2_cl}).
\end{example}

\begin{example}
Let $k=3,$ then $\left(k,12\right)=3,$ $D=1$. For $m\equiv-n\equiv1\mod4$
we have \[
\lambda_{m,3}=\frac{1}{c\left(0\right)}c\left(\frac{m-1}{4}\right),\,\mu_{n,3}=\frac{1}{c\left(0\right)}c\left(-\frac{m+1}{4}\right),\,\mbox{and}\,{\,\Lambda}_{3,R}=-R^{2}\left(1+R^{2}\right)^{2}.\]
For $m=n=3\equiv-1\mod4$ we get \begin{eqnarray*}
\mu_{3,3}^{2} & = & \Lambda_{3,R}\left(\lambda_{9,3}+\chi_{3}\left(3\right)\lambda_{1,3}\right),\end{eqnarray*}
or equivalently \[
c\left(-1\right)^{2}=\Lambda_{3,R}c\left(0\right)\left[c\left(2\right)-c\left(0\right)\right].\]
In Propositions \ref{pro:def_general_Tm} and \ref{pro:Th{m}_action}
one can compare any coefficient to compute eigenvalues $\lambda_{m}$
or $\mu_{m}$. Let $m\equiv1\mod4$ then using the first coefficient
gives the alternative expression \[
\lambda_{m,3}=\frac{1}{c\left(1\right)}\sum_{d|\left(m,\frac{5-m}{4}\right)}\chi_{3}\left(d\right)c\left(\frac{m}{d^{2}}+\frac{k\left(m-d^{2}\right)}{12d^{2}}\right)=\frac{1}{c\left(1\right)}\left[c\left(\frac{5m-1}{4}\right)+c\left(\frac{m-5}{20}\right)\right].\]
Changing variables to $l=\frac{m-1}{4}$ we get the following relation
for each positive integer $l$: \[
c\left(l\right)c\left(1\right)=c\left(0\right)\left[c\left(5l+1\right)+c\left(\frac{l-1}{5}\right)\right].\]
Table \ref{tab:relations} contains numerical values of some of the
above mentioned relations for $f\in\Mas{\Gamma,v,3.31105967012734}$
normalized with $c\left(1\right)=1$ and Fourier coefficients given
in Table \ref{tab:N1w3R3311}.

\end{example}

\begin{example}
Let $k=5$, then $\left(k,12\right)=1$ and $D=5$. For $m\equiv-n\equiv1\mod12$
with $\left(m,D\right)=\left(n,D\right)=1$ we have \begin{eqnarray*}
\lambda_{m,5} & = & \frac{1}{c\left(0\right)}c\left(5\frac{m-1}{12}\right),\,\mu_{n,3}=\frac{1}{c\left(0\right)}c\left(-5\frac{m+1}{12}\right),\,\mbox{and}\\
\Lambda_{5,R} & = & -R^{2}\left(1+R^{2}\right)^{2}\left(4+R^{2}\right)^{2}.\end{eqnarray*}
For $m=n=13$ we have $\lambda_{13}^{2}=\lambda_{169}+\chi_{5}\left(13\right)\lambda_{1}$
or equivalently \[
c\left(5\right)^{2}=c\left(0\right)\left[c\left(70\right)+c\left(0\right)\right].\]
For $n=11$ and $m\equiv-1\mod12$ we have \[
\mu_{11,5}\mu_{m,5}=\Lambda_{5,R}\left[\lambda_{11m}+\chi_{k}\left(11\right)\lambda_{\frac{m}{11}}\right],\]
where we also interpret $\lambda_{x}$ as $0$ unless $x$ is a positive
integer. Equivalently using change of coordinates to we have for all
$l\equiv0\mod5$ with $\left(\frac{12l}{5}-1,5\right)=1$: \[
c\left(-5\right)c\left(-l\right)=\Lambda_{5,R}c\left(0\right)\left[c\left(11l-5\right)-c\left(\frac{l-5}{11}\right)\right].\]
To illustrate relations involving the more complicated eigenvalues
possible for $D>1$ we consider $m=35\equiv-25\mod60$ and $n=11\equiv-1\mod12$.
We have \[
\mu_{35,5}\mu_{11,5}=\Lambda_{5,R}\lambda_{385},\]
and here $\mu_{35}=\frac{1}{c\left(0\right)}\left[c\left(-5\frac{35+1}{12}\right)+\chi_{5}\left(5\right)c\left(-5\frac{35+25}{12\cdot25}\right)\right]=\frac{1}{c\left(0\right)}\left[c\left(-15\right)+c\left(-1\right)\right]$
and $\lambda_{385}=\frac{1}{c\left(0\right)}\left[c\left(160\right)+c\left(6\right)\right]$.
The corresponding coefficient relation is thus \[
\left[c\left(-15\right)+c\left(-1\right)\right]c\left(-5\right)=-R^{2}\left(1+R^{2}\right)^{2}\left(4+R^{2}\right)^{2}c\left(0\right)\left[c\left(160\right)+c\left(6\right)\right].\]
As an example we consider $f\in\Mas{\Gamma,v,5,R}$ with $R=3.6624068669081$
and coefficients given by Table \ref{tab:N1w5R366}. Some numerical
examples of the above mentioned relation is found in Table \ref{tab:relations}.
\end{example}

\begin{example}
Let $k=6$, then $\left(k,12\right)=6$ and $D=\frac{k}{\left(k,12\right)}=1$.
This is a very special case since $k\equiv-k\mod12$ and thus $\overline{v}=v$.
For a positive odd integer $m$ both types of Hecke operators can
be defined and \begin{eqnarray*}
\lambda_{m,6} & = & \frac{1}{c\left(0\right)}c\left(\frac{m-1}{2}\right),\,\mu_{m,6}=\frac{1}{c\left(0\right)}c\left(-\frac{m+1}{2}\right),\mbox{ \mbox{and}}\\
\Lambda_{R,6} & = & \left(\frac{1}{4}+R^{2}\right)^{2}\left(\frac{9}{4}+R^{2}\right)^{2}\left(\frac{17}{4}+R^{2}\right)^{2}.\end{eqnarray*}
For $m=5,\, n=3$ we have $\lambda_{5,6}\lambda_{3,6}=\lambda_{15,6}$
and $\mu_{5,6}\mu_{3,6}=\Lambda_{6,R}\lambda_{15,6}$ or equivalently\begin{eqnarray*}
c\left(2\right)c\left(1\right) & = & c\left(0\right)c\left(7\right),\,\mbox{and}\\
c\left(-3\right)c\left(-2\right) & = & \Lambda_{6,R}c\left(0\right)c\left(7\right).\end{eqnarray*}
From this we see that \[
c\left(2\right)c\left(1\right)=c\left(-2\right)c\left(-3\right)\Lambda_{6,R}^{-1}=c\left(-2\right)c\left(-3\right).\]

\end{example}
For $m=1$ we have $\mu_{1,6}^{2}=\Lambda_{6,R}\lambda_{1,6}$ which
implies that \[
c\left(-1\right)^{2}=c\left(0\right)\Lambda_{6,R}.\]
For any positive integer $l$ set $m=1$ and $n=2l-1$ then $\mu_{1,6}\mu_{n,6}=\Lambda_{6,R}\lambda_{n,6}$
and hence \[
c\left(-1\right)c\left(-l\right)=c\left(0\right)\Lambda_{6,R}c\left(\frac{n-1}{2}\right)=c\left(-1\right)^{2}c\left(l-1\right),\]
which implies that \begin{equation}
c\left(-l\right)=c\left(-1\right)c\left(l-1\right).\label{eq:weight_k_prop_rel_c_-l_c_l-1}\end{equation}

For a specific example consider $f\in\Mas{\Gamma,v,6,R}$ with $R=3.70330780105981$,
normalized so that $c\left(0\right)=1$ and Fourier coefficients as
given in Table \ref{tab:N1w6R307}. See also in particular the third
column of this table where (\ref{eq:weight_k_prop_rel_c_-l_c_l-1})
is listed. 

\bibliographystyle{amsplain}
\bibliography{/home/fredrik/Documents/matematik/refs}

\begin{table}[b]
    \centering
    \caption[]{Fourier coefficients for 
$f \in  \MAS (\Gamma,v,2,2.95645894117486)$}
    \label{tab:N1w2R295}
    \begin{tabular}[t]{ld{15}d{15}d{10}}
$n$ &
\multicolumn{1}{c}{$c(n)$} &
\multicolumn{1}{c}{$c(-n)$} & 
\multicolumn{1}{c}{$|\frac{c(-1)c(-n)}{c(0)\Lambda_{2,R}(c(5n-1)}-1|$}  \\
 \hline\noalign{\smallskip}
  0&   1.230701624761 & \\
    1&   1.000000000000
 &  18.203610985364     & 1.98E-12\\
    2&  -1.471655144989
 &  13.721717970930     & 2.74E-10 \\
    3&   0.013247843393
 &  -3.974833500999     & 8.49E-10 \\
    4&   2.100340974458
 &  11.112158057949     & 3.49E-09 \\
    5&  -0.322346202386
 &  -14.371888713650    & 2.07E-09 \\
    6&  -1.246996804957
 &   14.791246405247    & 1.11E-09\\
    7&  -1.181658143432
 &   -6.282086679505    & 3.46E-08 \\
    8&  -0.418156991484
 &   -6.759700833281    & 2.92E-09 \\
    9&   2.510909902434
 &   10.100901280765    & 2.55E-09 \\
   10&   0.082126224604
 &    5.167895657723    & 3.27E-08 \\
   11&   1.165274980421
 &  -21.767613847686    & 1.02E-09\\
   12&  -0.247163744266
 &   13.245951772026    & 7.91E-09 \\
   13&  -2.102307863443
 &   11.149508424549    & 5.87E-09 \\
   14&  -0.727346882773
 &  -14.481437911536    & 2.05E-08 \\
   15&  -1.195785484887
 &   -1.591258609999    & 3.37E-07 \\
   16&   1.342209625858
 &    0.195951447099    & 3.40E-06 \\
   17&  -0.396934741247
 &    9.471610827564    & 5.46E-08 \\
   18&   1.423140585656
 &    7.748951323143    & 1.16E-07 \\
   19&   2.033391713377
 &  -21.084067140859    & 6.90E-08 \\
   20&   0.661999637495
 &   -3.229728578857    & 3.19E-07 \\
   21&  -1.566036852675
 &   12.863053232917    &  2.52E-07\\
   22&   0.010764537313
 &   12.313519090285    & 9.91E-07 \\
   23&  -0.162775774047
 &   -8.583584980476    & 1.98E-06 \\
   24&  -2.629883360741
 &  -16.408231323813    & 7.23E-09 \\
   25&   0.353626555761
 &   10.101583676954    & 2.92E-07 \\
\end{tabular}
\end{table}

\begin{table}[t]
    \centering
    \caption[]{Fourier coefficients for $f \in  \MAS (\Gamma,v,3,3.31105967012734)$}
    \label{tab:N1w3R3311}
    \begin{tabular}[t]{ld{15}d{15}d{12}}
$n$ &
\multicolumn{1}{c}{$c(n)$} &
\multicolumn{1}{c}{$c(-n)$} &
\multicolumn{1}{c}{$|\frac{c(n)c(1)}{c(0)\left(c(5n+1)+c\left(\frac{n-1}{5}\right)\right)}-1|$}  \\
 \hline\noalign{\smallskip}
     0&  1.531512741936 & \\
      1&  1.000000000000
 &   8.639530686828     & 4.32E-10 \\
      2& -1.574741762999
 &   5.624519754411     & 2.12E-10 \\
      3&  1.103109670076
 &  -2.889141766985     & 2.21E-10 \\
      4&  1.290948121024
 &   5.641174540769     & 2.03E-10 \\
      5& -2.022238297269
 &  -7.222396876842     & 9.59E-11 \\
      6& -8.785635544767
 &   5.995198221818     & 2.93E-10 \\
      7&  1.306633982800
 &  -2.438624530260     & 3.38E-10 \\
      8&  1.038761242319
 &  -5.703190028137     & 1.90E-11 \\
      9&  1.226237817255
 &   3.672525602758     & 4.88E-10 \\
     10& -2.479988645388
 &   6.222834185528     & 2.77E-10 \\
     11& -1.028226354383
 &  -1.026123492546     &  6.04E-11\\
     12&  2.149922429075
 &   7.860915325414     & 1.74E-09 \\
     13& -1.652728768231
 &   7.282463672083     & 3.97E-10 \\
     14&  2.596738601725
 &  -1.886462768742     & 7.84E-11 \\
     15&  2.661735950372
 &  -1.538694327170     & 2.95E-11 \\
     16&  7.202745626038
 &  -5.783279441143     &  1.78E-10\\
     17& -2.155511932252
 &   3.230287421009     & 4.81E-10 \\
     18& -1.550283931173
 &   9.405825626012     & 2.79E-10 \\
     19&  6.762558452563
 &  -4.956130358085     & 1.24E-10 \\
\end{tabular}
\end{table}

\begin{table}[b]
    \centering
    \caption[]{Fourier coefficients for 
$f \in  \MAS (\Gamma,v,5,3.6624068669081)$}
    \label{tab:N1w5R366}
    \begin{tabular}[t]{ld{19}d{19}}
     $n$ &
\multicolumn{1}{c}{$c(n)$} &
\multicolumn{1}{c}{$c(-n)$} \\
 \hline\noalign{\smallskip}
     0&  1.836234282360E+00  & \\
      1&  1.000000000000E+00
 &   1.887653985658E+03        \\
      2& -1.307103934737E+00
 &   1.044616926030E+03        \\
      3&  1.727292605687E+00
 &  -9.560107355852E+02       \\
      4& -8.244874356696E-01
 &   1.510043954372E+03        \\
      5& -1.357966916140E+00
 &  -1.555662019170E+03       \\
      6&  1.892760722402E+00
 &   3.969836744868E+02        \\
      7&  7.812520919140E-01
 &   1.363335177951E+03       \\
      8& -2.229697105219E+00
 &  -1.395993779667E+03       \\
      9& -1.107067514206E-01
 &  -1.001098116897E+03        \\
     10&  4.683894465654E-01
 &   2.325666345534E+03       \\
     11&  1.827070532237E+00
 &  -2.505722786009E+02       \\
     12&  2.758640358235E-01
 &  -1.190677542244E+03        \\
     13& -2.829618456564E+00
 &  -4.463696112166E+02        \\
     14&  1.075820859399E-01
 &   1.866263836513E+02       \\
     15&  7.874795758404E-02
 &   2.369156369887E+03        \\
     16&  1.583870281605E+00
 &  -8.471984859447E+02       \\
     17&  1.047442965873E+00
 &  -1.838202255925E+03       \\
     18& -7.395390192267E-01
 &  -7.508808946066E+01        \\
     19& -3.704603525272E-02
 &   4.924426686616E+02        \\
\end{tabular}
\end{table}

\begin{table}
    \centering
    \caption[]{Fourier coefficients for  
$f \in  \MAS (\Gamma,v,6,3.70330780105981)$}
    \label{tab:N1w6R307}
    \begin{tabular}[t]{ld{19}d{19}d{10}}
$n$ &
\multicolumn{1}{c}{$c(n)$} &
\multicolumn{1}{c}{$c(-n)$}  &
\multicolumn{1}{c}{$|\frac{c(-n)}{c(n-1)*c(-1)}-1|$}  \\
 \hline\noalign{\smallskip}
      0&  1.000000000000E+00  &   &        \\
      1&  5.384467700193E-01
 &   4.450801589309E+03       &          \\
      2& -5.994237858271E-01
 &   2.396519787896E+03       &         2.01E-08 \\
      3&  8.803713414082E-01
 &  -2.667916388933E+03       &         1.90E-08 \\
      4& -7.100750762253E-01
 &   3.918358240227E+03       &         1.91E-08 \\
      5& -2.098755988687E-01
 &  -3.160403332781E+03        &        1.74E-08 \\
      6&  1.016564356574E+00
 &  -9.341146729958E+02       &         2.57E-08 \\
      7& -3.227578035206E-01
 &   4.524526338100E+03       &         1.86E-08 \\
      8& -1.094555168655E+00
 &  -1.436530957159E+03        &        8.55E-09 \\
      9&  8.087258779256E-01
 &  -4.871647980231E+03        &        1.97E-08 \\
     10&  4.740331040530E-01
 &   3.599478507994E+03        &        2.37E-08 \\
     11&  2.179209189835E-02
 &   2.109827321226E+03       &         1.34E-08 \\
     12& -6.406911244172E-01
 &   9.699225130255E+01       &         2.68E-07 \\
     13& -9.207844021060E-01
 &  -2.851589119538E+03        &        1.57E-08 \\
     14&  1.219277100707E+00
 &  -4.098228742741E+03        &        1.52E-08 \\
     15&  5.223989524614E-01
 &   5.426760591533E+03        &        2.47E-08 \\
     16& -1.130068391137E-01
 &   2.325094105392E+03       &         7.54E-09 \\
     17& -5.277155225487E-01
 &  -5.029710458535E+02       &         5.31E-08 \\
     18& -8.607212197260E-01
 &  -2.348757106541E+03        &        8.55E-09 \\
     19&  5.473657942974E-01
 &  -3.830899448460E+03        &        1.98E-08 \\
     20& -2.046863808710E-01
 &   2.436216634091E+03        &        3.57E-08 \\
     21&  1.294360490006E+00
 &  -9.110184900836E+02       &         2.28E-08 \\
     22&  4.256358907948E-01
 &   5.760941791163E+03       &         1.13E-08 \\
     23& -1.407032267197E+00
 &   1.894420946282E+03       &         2.48E-08 \\
     24& -2.249462990120E-01
 &  -6.262421586507E+03       &         2.16E-08 \\
     25& -5.893596950834E-01
 &  -1.001191355161E+03        &        1.00E-08 \\
     26&  8.136637358514E-01
 &  -2.623123107600E+03       &         1.53E-08 \\
 
\end{tabular}
\end{table}

\begin{table}
    \centering
    \caption[]{Various Multiplicative Relations}
    \label{tab:relations}
\begin{tabular}[t]{cc} 

  \begin{tabular}[t]{lcd{6}}
    \multicolumn{3}{c}{$k=2,R=2.95645894117486$} \\
    \hline\noalign{\smallskip}
    $\left| c(1)c(2)-c(0)c(15)\right|$ & = & 5.9E-09  \\
    $\left| c(1)^{2}-c(0)\left( c(8)+c(0) \right) \right|$ & = & 3.7E-10 \\
    $\left|\frac{c\left(-1\right)c(-1996)}{c\left(0\right)\Lambda_{2,R}\left[c(9979)+c(399)\right]}-1\right|$
    & = & 4.7E-08\\
  \end{tabular}
  &
  \begin{tabular}[t]{lcd{6}}
    \multicolumn{3}{c}{$k=3,R=3.31105967012734$} \\
    \hline\noalign{\smallskip}
    $\left|\frac{c\left(-1\right)^{2}}{c\left(0\right)\Lambda_{6,R}\left[c\left(2\right)-c\left(0\right)\right]}-1\right|$
    & = & 7.8E-12 \\
    $\left|\frac{c\left(2\right)c\left(1\right)}{c\left(0\right)c\left(11\right)}-1\right|$
    & = & 2.1E-10 \\
    $\left|\frac{c\left(3\right)c\left(1\right)}{c\left(0\right)c\left(16\right)}-1\right|$
    & = & 2.2E-10 \\
  \end{tabular}
  \\
  \begin{tabular}[t]{lcd{6}}
    \multicolumn{3}{c}{$k=5,R=3.6624068669081$} \\
    \hline\noalign{\smallskip}
    $\left| c(5)^{2} - c(0)\left( c(70)+c(0) \right) \right|$ & = & 5.4E-08 \\
    $\left| \frac{\left[c\left(-15\right)+c\left(-1\right)\right]c\left(-5\right)}{\Lambda_{5,R}c\left(0\right)\left[c\left(160\right)+c\left(6\right)\right]}-1\right|$
    & = & 7.7E-07 \\
    $\left| \frac{c(-5)c(-5)}{\Lambda_{5,R}\left(c(50)-c(0)\right)}-1\right|$
    & = & 2.1E-08  \\
  \end{tabular}
  &
  \begin{tabular}[t]{lcd{6}}
    \multicolumn{3}{c}{$k=6,R=3.70330780105981$} \\
    \hline\noalign{\smallskip}
    $\left|c\left(2\right)c\left(1\right)-c\left(0\right)c\left(7\right)\right|$
    & = & 2.1E-09 \\
    $\left|c\left(-1\right)^{2}/\Lambda_{6,R}-1\right|$ & = & 3.6E-08 \\
    $\left|\frac{c\left(-139\right)}{c\left(138\right)c\left(-1\right)}-1\right|$
    & = & 2.3E-07 \\
  \end{tabular}
  \\
  \end{tabular}
\end{table}
\end{document}